\documentclass{amsart}
\usepackage{graphicx}
\newtheorem{thm}{Theorem}[section]

\newtheorem{lem}[thm]{Lemma}

\theoremstyle{definition}
\newtheorem{defn}[thm]{Definition}
\theoremstyle{example}
\newtheorem{ex}[thm]{Example}

\numberwithin{equation}{section}
\begin{document}

\title[Some results in Generalized \v{S}erstnev spaces]
{Some results in Generalized \v{S}erstnev spaces}%
\author[Alimohammady and Saadati]{M. Alimohammady and R. Saadati }%
\address{Dept. of Math., Azad University, Nur, Iran}
\email{amohsen@umz.ac.ir}
\address{Dept. of Math., Azad University, Amol, P.O.Box 678, Iran}%

\email{rsaadati@eml.cc}%
\thanks{}%

%\date{}%
%\dedicatory{}%
%\commby{}%
% ----------------------------------------------------------------
\begin{abstract}
In this paper, we show that  D-compactness in Generalized
\v{S}erstnev spaces implies D-boundedness and as in the classical
case, a D-bounded and closed subset of a
 characteristic Generalized \v{S}erstnev is not
D-compact in general. Finally, in the finite dimensional
Generalized \v{S}erstnev spaces a subset is D-compact if and only
if it is D-bounded and closed.
\end{abstract}
\maketitle
% ----------------------------------------------------------------
\section{Introduction and Preliminaries}
Probabilistic normed spaces (PN spaces henceforth) were introduced
by \v{S}erstnev in 1963 \cite{ser}. In the sequel, we adopt the
new definition of Generalized \v{S}erstnev  PN spaces given in the
paper by  Lafuerza-Guill\'{e}n and Rodr\'{i}guez \cite{gu3}. The
notations and concepts used are those of \cite{al1,al2,gu1,gu3}
and \cite{sch1}.

In the sequel, the space of probability distribution functions is
denoted by $ \Delta^{+} =\{ F:\mathbb{R}\longrightarrow [0,1]: F$
 is left-continuous, non-decreasing,  $F(0)=0$ and
$F(+\infty)=1\}$ and  $D^{+} \subseteq \Delta^{+}$ is defined as
follows
$$ D^{+}=\{F\in
 \Delta^{+}:l^{-}F(+\infty)=1\}.$$ The space $\Delta^{+}$ is
partially ordered by the usual point-wise ordering of functions
i.e., $F\leq G$ if and only if $F(x)\leq G(x)$ for all $x$ in
$\mathbb{R}$. The maximal element for $\Delta^{+}$ in this order
is $\varepsilon_{0}$ , a distributive defined by
$$
\varepsilon_{0}=\begin{cases} \begin{array}{ccc}
  0 & if & x\leq 0 \\
1 & if & x>0 .\\
\end{array}
\end{cases}$$
 A triangle function is a binary
operation on $\Delta^{+}$, namely a function
$\tau:\Delta^{+}\times \Delta^{+}\longrightarrow \Delta^{+}$ that
is associative, commutative, nondecreasing and which has
$\varepsilon_{0}$ as unit, viz. For all $F,G,H\in\Delta^{+}$, we
have
\begin{eqnarray*} \tau(\tau(F,G),H)&=&\tau(F,\tau(G,H)),\\
\tau(F,G)&=&\tau(G,F),\\
F\leq G &\Longrightarrow & \tau(F,H)\leq\tau(G,H),\\
\tau(F,\varepsilon_{0})&=&F.\end{eqnarray*} Continuity of a
triangle functions means continuity with respect to the topology
of weak convergence in $\Delta^{+}$.

Typical continuous triangle functions are
$\tau_{T}(F,G)(x)=\sup_{s+t=x}T(F(s),G(t))$ and
$\tau_{T^{*}}(F,G)=\inf_{s+t=x}T^{*}(F(s),G(t)).$

Here, $T$ is a continuous t-norm, i.e. a continuous binary
operation on $[0,1]$ that is commutative, associative,
nondecreasing in each variable and has 1 as identity, and $T^{*}$
is a continuous t-conorm, namely a continuous binary operation on
$[0,1]$ which is related to the continuous t-norm $T$ through $
T^{*}(x,y)=1-T(1-x,1-y).$

\begin{defn}A probabilistic normed (briefly PN) space is a
quadruple $(V,\nu,\tau,\tau^{*})$, where $V$ is a real vector
space, $\tau$ and $\tau^{*}$ are continuous triangle functions,
and $\nu$ is a mapping from $V$ into $\Delta^{+}$ such that, for
all $p,q$ in $V$, the following conditions hold:

\,  (N1) $\nu_{p}=\varepsilon_{0}$ if and only if, $p=\theta$,
where $\theta$ is the null vector in $V$;

\,  (N2) $\nu_{-p}=\nu_{p},$ for each $p\in V$;

\,  (N3) $\nu_{p+q}\geq \tau(\nu_{p},\nu_{q}),$ for all $p,q\in
V$;

\,  (N4) $\nu_{p}\leq \tau^{*}(\nu_{\alpha p},\nu_{(1-\alpha)p})$,
for all $\alpha$ in $[0,1]$.\end{defn}If the inequality (N4) is
replaced by the equality $$\nu_{p}=\tau_{M}(\nu_{\alpha
p},\nu_{(1-\alpha)p}),$$ then the PN space is called
\textit{\v{S}erstnev } space and, as a consequence, a condition
stronger than (N2) holds, namely $$ \nu_{\lambda
p}(x)=\nu_{p}(\frac{x}{|\lambda|}),$$ for all $p\in V$ , $\lambda
\neq 0$ and $x\in \mathbb{R}$.

Following \cite{al0,gu3}, for $0 < b \leq +\infty$, let $M_b$ be
the set of \emph{m-transforms}  consist on all continuous and
strictly increasing functions from $[0, b]$ onto $[0,+\infty]$.
More generally, let $\widetilde{M}$ be the set of non decreasing
left-continuous functions $\phi : [0,+\infty] \longrightarrow
[0,+\infty]$, $\phi(0) = 0$, $\phi(+\infty) = +\infty$ and
$\phi(x)> 0$, for $x> 0$. Then $M_b\subseteq\widetilde{M}$ once
$m$ is extended to $[0,+\infty]$ by $m(x) = +\infty$ for all $x
\geq b$. Note that a function $\phi \in \widetilde{M}$ is
bijective if and only if $\phi\in M_{+\infty}$. Sometimes, the
probabilistic norms $\nu$ and $\nu'$ of two given  PN spaces
satisfy $\nu'=\nu\phi$ for some $\phi\in M_{+\infty}$, non
necessarily bijective. Let $\hat{\phi}$ be the (unique)
quasi-inverse of $\phi$ which is left-continuous. Recall from
\cite{sch1} page 49, that $\hat{\phi}$ is defined by
$\hat{\phi}(0) = 0$, $\hat{\phi}(+\infty) = +\infty$ and
$\hat{\phi}(t) = \sup\{u : \phi(u) < t\}$ for all $0 < t <
+\infty$. It follows that $\hat{\phi}(\phi(x)) \leq x$ and
$\phi(\hat{\phi}(y)) \leq y $ for all $x$ and $y$.
\begin{defn}\cite{gu3}.  A quadruple $(V, \nu,\tau,\tau^*)$ satisfying the $\phi$-\v{S}erstnev condition
\begin{eqnarray*} \nu_{\lambda p}(x) = \nu_{p}
(\hat{\phi}(\frac{\phi(x)}{|\lambda|})),\end{eqnarray*} for all
$x\in \mathbb{R}^+$ , $p\in V$ and $\lambda\in \mathbb{R}- \{0\}$
 is called a $\phi$-\v{S}erstnev PN
space (\emph{Generalized} \v{S}erstnev spaces).\end{defn}
\begin{lem} If $|\alpha|\leq |\beta|$, then $\nu_{\beta p}\leq
\nu_{\alpha p}$ for every $p$ in $V$.\end{lem}
\begin{defn} Let $(V,\nu,\tau,\tau^{*})$ be a PN space. For each $p$ in
$V$ and $\lambda>0$, the strong $\lambda-neighborhood$ of $p$ is
the set
\begin{eqnarray*} N_{p}(\lambda)&=&\{q\in V:\nu_{p-q}(\lambda)>1-\lambda\},
\end{eqnarray*}
and the strong neighborhood system for $V$ is the union
$\bigcup_{p\in V}\mathcal{N}_{p}$ where
$\mathcal{N}_p=\{N_p(\lambda):\lambda>0\}$.
\end{defn}

The strong neighborhood system for $V$ determines a Hausdorff
topology for $V$.

\begin{defn} Let $(V,\nu,\tau,\tau^{*})$ be a PN space, a sequence
$\{p_{n}\}_{n}$ in $V$ is said to be strongly convergent to $p$ in
$V$ if for each $\lambda>0$, there exists a positive integer $N$
such that $p_{n}\in N_{p}(\lambda )$, for $n\geq N$. Also the
sequence $\{p_{n}\}_{n}$ in $V$ is called strongly Cauchy
 sequence if for every $\lambda >0$ there is a positive
 integer $N$
 such that $\nu_{p_{n}-p_{m}}(\lambda )>1-\lambda $, whenever $m,n>N$.
 A PN space $(V,\nu,\tau,\tau^{*})$ is said to be
 strongly complete in the strong topology if and only if every strongly Cauchy sequence in $V$ is
 strongly convergent to a point in $V$.
\end{defn}
\begin{defn}\cite{gu1}. Let $(V,\nu,\tau,\tau^{*})$ be a PN space and $A$ be
 the nonempty subset of $V$. The probabilistic radius of $A$ is
 the function $R_{A}$ defined on $\mathbb{R}^{+}$ by
 $$R_{A}(x)=\begin{cases}\begin{array}{ccc}
   l^{-}inf\{\nu_{p}(x):p\in A\} &  if& x\in[0,+\infty) \\
   1 &  if& x=+\infty. \\
 \end{array}
\end{cases}$$
\end{defn}

\begin{defn}\cite{gu1}. A nonempty set $A$ in a PN space
$(V,\nu,\tau,\tau^{*})$is said to be:

\, (a) certainly bounded, if $R_{A}(x_{0})=1$ for some $x_{0}\in
(0,+\infty )$;

\, (b) perhaps bounded, if one has $R_{A}(x)<1$, for every $x\in
(0,+\infty )$ and  $l^{-}R_{A}(+\infty )=1$;

\, (c) perhaps unbounded, if $R_{A}(x_{0})>0$ for some $x_{0}\in
(0,+\infty )$ and $l^{-}R_{A}(+\infty )\in(0,1)$;

\, (d) certainly unbounded, if $l^{-}R_{A}(+\infty )=0$ i.e., if
$R_{A}=\varepsilon_{\infty}$.\\
Moreover, $A$ is said to be distributionally bounded, or simply
D-bounded if either (a) or (b) holds, i.e. $R_{A} \in D^{+}.$ If
$R_{A}\in \Delta^{+} - D^{+}$ , $A$ is called
D-unbounded.\end{defn}

\begin{thm}\cite{gu1}. A subset $A$ in the PN space $(V,\nu,\tau,\tau^{*})$
is D-bounded if and only if there exists a d.f. $G\in D^{+}$ such
that $\nu_{p}\geq G$ for every $p\in A$.\end{thm}

\begin{defn}\cite{gu1}. A subset $A$ of TVS $V$ is said to be topologically
bounded if for every sequence $\{\alpha_{n}\}$ of real numbers
that converges to zero as $n\longrightarrow +\infty$ and for every
$\{p_{n}\}$ of elements of $A$, one has
$\alpha_{n}p_{n}\longrightarrow \theta$, in the strong
topology.\end{defn}
\begin{thm}\cite{al2}.Suppose $(V,\nu,\tau,\tau^{*})$ is a PN space, when it is
endowed with the strong topology induced by the probabilistic norm
$\nu$. Then it is a topological vector space if and only if for
every $p\in V$ the map from $\mathbb{R}$ into $V$ defined by
$$\alpha\longmapsto\alpha p$$ is continuous.\end{thm}
 The PN space $(V,\nu,\tau,\tau^{*})$ is called
\textit{characteristic} whenever $\nu(V)\subseteq D^+$.
\begin{thm}\cite{gu3}. Let $\phi\in \widetilde{M}$ such that
$\lim_{x\longrightarrow\infty}\hat{\phi}(x)=\infty$. Then a
 $\phi$-\v{S}erstnev  PN space $(V,\nu,\tau,\tau^*)$ is a TVS if and only if it is
characteristic.\end{thm}

\begin{lem}\cite{gu2,gu3}. Let $\phi\in \widetilde{M}$ such that
$\lim_{x\longrightarrow\infty}\hat{\phi}(x)=\infty$.
 Let $(V,\nu,\tau,\tau^*)$ be a characteristic $\phi$-\v{S}erstnev  PN
 space. Then for a subset $A$ of $V$ the following are equivalent

 \, (a) For every $n\in \mathbb{N}$, there is a $k\in\mathbb{N}$ such that $A\subset kN_{\theta}(1/n)$.

\, (b) $A$ is D-bounded.

\, (c) $A$  is topologically bounded .\end{lem}
\section{D-bounded and D-compact sets in $\phi$-\v{S}erstnev spaces}

\begin{thm} Let $\phi\in \widetilde{M}$ such that
$\lim_{x\longrightarrow\infty}\hat{\phi}(x)=\infty$. Then in a
characteristic $\phi$-\v{S}erstnev  PN space $(V,\nu,\tau,\tau^*)$
if $ p_{m}\longrightarrow p$ in $V$ and
$A=\{p_{m}:m\in\mathbb{N}\}$,
 then $A$ is a D-bounded subset of $V$.\end{thm}
 \textbf{Proof.} Let
$p_{m}\longrightarrow p_0$ and $\alpha_{m}\longrightarrow 0$. Then
there exists  $m_{0}\in \mathbb{N}$ such that for every $m\geq
m_{0}$,
  $0<\alpha_{m}<1$, then
\begin{eqnarray*} \nu_{\alpha_{m}p_{m}}&\geq
&\tau(\nu_{\alpha_{m}(p_{m}-p_{0})},\nu_{\alpha_{m}p_{0}})\\
&>& \tau(\nu_{p_{m}-p_{0}},\nu_{\alpha_{m}p_{0}}) \\
&\longrightarrow& \tau(\varepsilon_{0},\varepsilon_{0})\\
&=&\varepsilon_{0},\end{eqnarray*} as $m$ tend to infinity. \qed
\begin{ex} The 4-tuple
$(\mathbb{R},\nu,\tau_{\pi},\tau_{M})$, where
$\nu:\mathbb{R}\longrightarrow \Delta^{+}$ is defined by
$$ \nu_{p}(x)=\begin{cases}\begin{array}{ccc}
  0 & if & x=0 \\
  \frac{ax}{x+|p|} & if & 0<x<+\infty \\
  1 & if & x=+\infty \\
\end{array}\end{cases},$$
$a\in(0,1)$, $\nu_{0}=\varepsilon_{0}$ and  $\phi(x)=x$, is a
$\phi$-\v{S}erstnev space (see \cite{gu3}). The sequence
$\{\frac{1}{m}\}$ is convergent but $A=\{\frac{1}{m}: m\in
\mathbb{N}\}$ is not D-bounded set. The only D-bounded set in this
space is $\{0\}$.\end{ex}
\begin{defn}The characteristic  $\phi$-\v{S}erstnev  PN space $(V,\nu,\tau,\tau^*)$ is said to be
distributionally compact (simply D-compact) if every sequence
$\{p_{m}\}_{m}$ in $V$ has a convergent subsequence
$\{p_{m_{k}}\}$. A subset $A$ of a characteristic
$\phi$-\v{S}erstnev  PN space $(V,\nu,\tau,\tau^*)$ is said to be
D-compact if every sequence $\{p_{m}\}$ in $A$ has a subsequence
$\{p_{m_{k}}\}$ convergent to a element $p\in A$.\end{defn}

\begin{lem} A D-compact subset of a characteristic  $\phi$-\v{S}erstnev  PN space $(V,\nu,\tau,\tau^*)$,
is D-bounded and closed.\end{lem} \textbf{Proof.} Suppose that
$A\subseteq V$ is D-compact. From Lemma 1.12  it is enough show
that $A$ is topologically bounded. On the contrary there is a
sequence $\{p_{m}\}\subseteq A$ and a real sequence
$\alpha_{m}\longrightarrow 0$ such that $\{\alpha_{m} p_{m}\}$
doesn't tend to the origin in $V$. Then there is an infinite set
$J\subseteq \mathbb{N}$ such that the sequence
$\{\alpha_{m}p_{m}\}_{m\in J}$ lies in the complement of a
neighborhood of the origin. Now $\{p_{m}\}$ is a subset of
D-compact set $A$, so it has a convergent subsequence
$\{p_m\}_{m\in J'}$. From Lemmas 1.12 and 2.1 $\{p_m\}_{m\in J'}$
is topologically bounded and so $\{\alpha_{m}p_{m}\}_{m\in J'}$
tends to origin  which is a
contradiction. The closedness of $A$ is trivial .\qed \\

As in the classical case, a D-bounded and closed subset of a
 characteristic  $\phi$-\v{S}erstnev is not
D-compact in general, as one can see from the next example.

\begin{ex} We consider quadruple $(\mathbb{Q},\nu,\tau_{\pi},\tau_{M})$, where
$\pi(x,y)=xy$, for every $x,y\in[0,1]$ and probabilistic norm
$\nu_{p}(t)=\frac{t}{t+|p|}$. It is straightforward to check that
$(\mathbb{Q},\nu,\tau_{\pi},\tau_{M})$ is a characteristic
$\phi$-\v{S}erstnev  PN space. In this space, convergence of a
sequence is equivalent to its convergence in $\mathbb{R}$. We
consider the subset $A=[a,b]\cap \mathbb{Q}$ which $a,b\in
\mathbb{R}-\mathbb{Q}$. Since
$R_{A}(t)=\frac{t}{t+\max\{|a|,|b|\}}$, then $A$ is D-bounded set
and since $A$ is closed in $\mathbb{Q}$ classically, and so is
closed in $(\mathbb{Q},\nu,\tau_{\pi},\tau_{M})$. We know $A$ is
not classically compact in $\mathbb{Q}$, i.e. there exists a
sequence in $\mathbb{Q}$ with no convergent subsequence in
classical sense and so in $(\mathbb{Q},\nu,\tau_{\pi},\tau_{M})$.
Hence $A$ is not D-compact.\end{ex}
\begin{thm}\cite{sa}.Consider a finite dimensional characteristic  $\phi$-\v{S}erstnev  PN space $(V,\nu,\tau,\tau^*)$
on real field $(\mathbb{R},\nu',\tau',\tau'^{*})$. Every subset
$A$ of $V$ is D-compact if and only if $A$ is D-bounded and
closed.\end{thm}

% ----------------------------------------------------------------
\bibliographystyle{amsplain}
\bibliography{}
\end{document}